\documentstyle[12pt]{article}
\input{amssym}

\newcommand{\RR}{\Bbb{R}}
\newcommand{\SS}{\Bbb{S}}
\begin{document}
\title{ Fuzzy  Lie Groups}
\author{M. Nadjafikhah\thanks{e-mail: m\_nadjafikhah@iust.ac.ir} \and R. Bakhshandeh
Chamazkoti\thanks{e-mail: r\_bakhshandeh@iust.ac.ir}}
\date{\it Department of Mathematics, Iran University of Science and
Technology, Narmak, Tehran, IRAN.}
\maketitle
\renewcommand{\sectionmark}[1]{}

\begin{abstract}
In this paper, we have tried to apply the concepts of fuzzy sets
to Lie groups and its relative concepts. First, we define a ${\cal
C}^1$ fuzzy submanifold after reviewing ${\cal C}^1-$fuzzy
manifold definition. In main section, we defined the Lie group and
some its relative concepts such as fuzzy transformation group,
fuzzy $G$-invariant. Our aim is to prepare the suitable conditions
for defining fuzzy differential invariant by constructing suitable
objects.

\end{abstract}

\section{Introduction}

The notion of a fuzzy Lie group is depend on the basic concepts in
fuzzy topology \cite{[4]}, ${\cal C}^1-$fuzzy manifold and fuzzy
differentiable function between two ${\cal C}^1-$fuzzy manifolds
\cite{[12]}.

In this paper, we introduce the  fuzzy sets with two notations
which they are equivalent.

First, let $X$ be a non-empty set of points and $I=[0,1]$; then
$I^X$ will ne denoted the set of all functions
$\mu:X\longrightarrow I$. A member of $I^X$ is called {\it a fuzzy
subset of $X$}. If $x\in X$ and $p\in[0,1]$, by the fuzzy point
$x_p$, we mean the fuzzy subset of $X$ which takes the value $p$
at the point $x$ and 0 elsewhere. Let $\mu$ be a fuzzy set of $X$
and $t\in[0,1]$, the set $\{x\in X: \mu(x)\geq t\}$ is called {\it
a level subset of $\mu$} and is symbolized by $\mu_t$. An element
$a\in X$ is called {\it a normal element of $\lambda$ with respect
to $\mu$}, if $\lambda(a)\geq\mu(y)$ for all $y\in X$,
\cite{[13]}.

Second, a fuzzy set $A$ in $X$ is characterized by a membership
(characteristic) function $\mu_A:X\longrightarrow I$, which
associates with each points in $X$, a real number in interval
$[0,1]$, with the value $\mu_A(x)$ at $x$ representing the grade
of membership of $x$ in $A$. We have:
\begin{itemize}
\item[i)] if $A\subseteq B$, then $\mu_A\leq\mu_B$;
\item[ii)] if $C=A\cup B$, then $\mu_C(x)=\max\{\mu_A(x),\mu_B(x):x\in
X\}$; and
\item[iii)] if $D=A\cap B$, then $\mu_D(x)=\min\{\mu_A(x),\mu_B(x):x\in
X\}$.
\end{itemize}
more generally, for a family of fuzzy sets, $A=\{A_i: i\in J\}$,
the union, $C=\cup_{i\in J}A_i$, and intersection, $D=\cap_{i\in
J}A_i$, are defied by
\begin{eqnarray*}
\mu_C(x)&=&\sup\{\mu_{A_i}(x):i\in J,x\in X\},\\
\mu_D(x)&=&\inf\{\mu_{A_i}(x):i\in J,x\in X\},
\end{eqnarray*}
\paragraph{Convention.}
We will consider $I=[0,1]\subset\RR$ and $J$  an indexing set.

\section{Preliminaries}
In this section, we define some fuzzy concepts        use.
\paragraph{Definition 2.1.}
If $\lambda\in I^X$ , $\mu\in I^Y$ then $\lambda\times\mu\in
I^{X\times Y}$ is defined by
$$(\lambda\times\mu)(x,y):=\min\{\lambda(x),\mu(y)\},\hspace{1cm} \forall\,(x,y)\in X\times Y.$$

Let $F:\lambda\longrightarrow\mu$ be a fuzzy function. If
$A\leq\lambda$ , $B\leq\mu$, then $F(A)$, $F^{-1}(B)$ are defined
by
\begin{itemize}
\item[] $(F(A))(y):=\sup\Big\{\min\{F(x,y),A(x)\}\,:\,x\in X\Big\}$
for all $y\in Y$, and
\item[] $(F^{-1}(B))(x):=\sup\Big\{\min\{F(x,y),B(x)\}\,:\,y\in
Y\Big\} $ for all $x\in X$.
\end{itemize}

\paragraph{Definition 2.2.}
A fuzzy  function $F:\lambda\longrightarrow\mu$ is said to be
\begin{itemize}
\item[i)] {\it injective} if $F(x_1)=F(x_2)$
 then $ x_1=x_2$ ;
\item[ii)] {\it surjective} if for all $y\in Y$ with $\mu(y)\neq0$, there exists $x\in X$ such that $F(x,y)=\lambda(x)$;
\item[iii)] {\it bijective} if $F$ is both injective and surjective.
\end{itemize}
\paragraph{Definition 2.3. (Chakrabarty and Ahsanullah
\cite{[2]})}
Let $\mu$ be a fuzzy subset of $X$. A collection $\tau$ of fuzzy
subsets of $\mu$ satisfying:
\begin{itemize}
\item[i)] $t\bigcap\mu\in\tau$ for all $t\in I$,
\item[ii)] if $\mu_i\in\tau$ for all $i\in J$ then
$\bigcup\{\mu_i:i\in J\}\in\tau$,
\item[iii)] if $\mu,\nu\in\tau$ then $\mu\bigcap\nu\in\tau$.
\end{itemize}
$\tau$ is called {\it a fuzzy topology on} $\mu$. The pair
$(\mu,\tau)$ is called {\it a fuzzy topological space}. Members of
$\tau$ will be called {\it fuzzy open sets} and their complements
with respect to $\mu$ are called { \it fuzzy closed sets of
$(\mu,\tau)$}. ${\cal B}\subset\tau$ is called {\it an open base
of $\tau$} if every member of $\tau$ can be expressed as union of
some members of ${\cal B}$.

Let $(\lambda,\tau)$ and $(\mu,\tau^{\prime})$ are two fuzzy
topological spaces, the collection
$${\cal B}=\{\gamma\times\eta:
\gamma\in\tau,\eta\in\tau^{\prime}\}$$ form an open base of a
fuzzy topology in $\lambda\times\mu$. The fuzzy topology in
$\lambda\times\mu$, induced by ${\cal B}$ is called {\it the
product fuzzy topology} of $\tau$ and $\tau^{\prime}$ and is
denote by $\tau\times\tau^{\prime}$. The fuzzy topological space
$(\lambda\times\mu,\tau\times\tau^{\prime})$ is called {\it the
product of the fuzzy topological spaces $(\lambda,\tau)$ and
$(\mu,\tau^{\prime})$.}

A fuzzy topological space is called {\it a fuzzy $T_1-$space}, if
every fuzzy point is a closed fuzzy set. (\cite{[2]} and
\cite{[12]})
\paragraph{Definition 2.4.}
$(\lambda,\tau)$ is said to be a fuzzy Hausdorff space if
$x_p,y_p\in\lambda~~(x\neq y)$, there are exist
$\mu,\upsilon\in\tau$ such that $x_p\in\mu,y_q\in\upsilon$ and
$\mu\cap\upsilon=0$.
\paragraph{Definition 2.5.}
A fuzzy proper function
$F:(\lambda,\tau)\longrightarrow(\mu,\tau^{\prime})$ is said to be
\begin{itemize}
\item[i)] {\it fuzzy continuous} if
$F^{-1}(\nu)\in\tau$ for all $\nu\in\tau^{\prime}$,
\item[ii)] {\it fuzzy open} if
$F(\delta)\in\tau^{\prime}$ for all $\delta\in\tau$,
\item[iii)] {\it fuzzy homomorphism} if $F$ be bijective, fuzzy continuous and
open.
\end{itemize}
\paragraph{Definition 2.6.}
Let $G$ be a group. $\mu\in I^G$ is said to be {\it a fuzzy
subgroup of $G$} if $\forall x,y\in G$
\begin{itemize}
\item[i)] $\mu(xy)\geq\min\{\mu(x),\mu(y)\}$,
\item[ii)] $\mu(x^{-1})=\mu(x)$.
\end{itemize}
\paragraph{Definition 2.7. (Das \cite{[11]}).}
A fuzzy topology $\tau$ on a group $G$ is said to be {\it
compatible} if the mappings
\begin{eqnarray}
&m:(G\times G,\tau\times\tau)\longrightarrow (G,\tau)
\hspace{1cm} (x,y)\mapsto xy & \label{eq:1}\\
&i: (G,\tau)\longrightarrow (G,\tau) \hspace{1cm} x\mapsto x^{-1}&
\label{eq:2}
\end{eqnarray}
are fuzzy continuous. A group $G$ equipped  with a compatible
fuzzy topology $\tau$ on $G$ is called {\it a fuzzy topological
group.}
\paragraph{Definition 2.8.}
{\it A fuzzy topological vector space}, is a vector space $E$ over
the field $F$ of real or complex numbers, if $E$ is equipped with
a fuzzy topology $\tau$ and $F$ equipped with the usual topology
${\cal F}$, such that following two mappings are fuzzy continuous:
\begin{eqnarray*}
&(E,\tau)\times(E,\tau)\longrightarrow (E,\tau)\hspace{1cm}
(x,y)\mapsto x+y& \\
&(F,{\cal F})\times(E,\tau)\longrightarrow (E,\tau)\hspace{1cm}
(\alpha,x)\mapsto \alpha x.&
\end{eqnarray*}
\paragraph{Definition 2.9.}
Let $E,F$ be two fuzzy topological vector space, the mapping
$\phi:E\longrightarrow F$ is said to be {\it tangent to $0$} if
given a neighborhood $W$ of $0_{\delta}$, $0<\delta\leq1$, in $F$
there exists a neighborhood $V$ of $0_{\epsilon}$,
$0<\epsilon<\delta$, in $E$ such that
$$\phi[tV]\subset o(t)W,$$
for some function $o(t)$.
\paragraph{Definition 2.10. (Ferraro and Foster \cite{[12]}).}
Let $E,F$ be two fuzzy topological vector space, each endowed with
a $T_1-$fuzzy topology. Let $f:E\longrightarrow F$ be a fuzzy
continuous mapping. The $f$ is said to be {\it fuzzy
differentiable at a point $x\in E$} if there exists a linear fuzzy
continuous mapping $u:E\longrightarrow F$ (It is denoted by
$u\in{\cal L}(E,F).$) such that
$$f(x+y)=f(x)+u(y)+\phi(y),~~~~~y\in E,$$
where $\phi$ is tangent 0. The mapping $u$ is called {\it the
fuzzy derivative of $f$ in $x$} that is denoted by
$f^{\prime}(x)$; $f^{\prime}(x)\in{\cal L}(E,F)$. The mapping $f$
is {\it fuzzy differentiable} if it is fuzzy differentiable at
every point of $E$.
\paragraph{Definition 2.11.}
A fuzzy proper function $F:\lambda\longrightarrow\mu$ is said to
be {\it a fuzzy homomorphism} if $F(x,y)=\lambda(x)$,
$F(z,w)=\lambda(z)$ imply
 $$F(xz,yw)=\lambda(xz) \hspace{1cm} \mbox{for all}\;\; x,y,z,w\in G.$$
\paragraph{Definition 2.12.}
Let $E,F$ be fuzzy topological vector spaces. A bijection
$f:E\longrightarrow F$ is said to be a {\it ${\cal C}^1-$fuzzy
diffeomorphism} if it and its inverse $f^{-1}$ are fuzzy
differentiable, and $f^{\prime}$ and $(f^{-1})^{\prime}$ are fuzzy
continuous.
\paragraph{Definition 2.13. (Kim and Lee \cite{[5]})}
Let $X$ be a vector space over a field $F$ with Lie bracket
$X\times X\ni(x,y)\mapsto[x,y]\in X$. $X$ is called {\it a Lie
algebra over $F$} if
\begin{itemize}
\item[i)] The Lie bracket is bilinear,
\item[ii)] $[x,y]=-[y,x]$ for all $x,y\in X$,
\item[iii)] (Jacobi identity ) $[x,[y,z]]+[y,[z,x]]+[z,[x,y]]=0$ for all $x,y,z\in
X$.
\end{itemize}
 Now  $\mu\in I^X$  is called {\it a
fuzzy Lie subalgebra of $X$} if, for all $\alpha\in F, x,y\in X$,
the following requirements are met
\begin{itemize}
\item[i)] $\mu(x+y)\geq\min\{\mu(x),\mu(y)\}$,
\item[ii)] $\mu(\alpha x)\geq\mu(x)$,
\item[iii)] $\mu([x,y])\geq\min\{\mu(x),\mu(y)\}$.
\end{itemize}
If the condition (iii) is replaced by
$$\mu([x,y])\geq\max(\mu(x),\mu(y)),$$
then $\mu$ is called {\it a fuzzy Lie ideal of $X$}.
\paragraph{Example 2.14.}
Let $X={\RR}^3$ and $[x,y]=x\times y$, where $\times$ is cross
product, for all $x,y\in X$. Then $X$ is a Lie algebra over a
field ${\RR}$. Define $\mu:{\RR}^3\longrightarrow I$ by
$$f(x)=\left\{ \begin{array}{lcl}
1 && \mbox{if}\;x=y=z=0, \\ 1/4 && \mbox{if}\;x=y=0\;\mbox{and}\;
z\neq0 \\ 0 && \mbox{otherwise}.
\end{array}\right.$$
then $\mu$ is a fuzzy subalgebra of $X$. But $\mu$ is not a fuzzy
Lie ideal of $X$; because
$$\mu([(0,0,1),(1,1,1)])=\mu((0,0,1)\times(1,1,1))=\mu(-1,1,0)=0,$$
while
$$\max\{\mu(0,0,1),\mu(1,1,1)\}=\max\Big\{{1\over4},0\Big\}=\frac{1}{4}.$$
\section{ ${\cal C}^1-$Fuzzy Submanifold}
We Start this section by defining a ${\cal C}^1-$fuzzy manifold
and some examples:
\paragraph{Definition 3.1. (Ferraro and Foster \cite{[12]})}
Let $X$ be a set. A {\it ${\cal C}^1-$fuzzy atlas on $X$} is a
collection of pairs $\{(A_j,\phi_j)\}_{j\in J}$, which satisfies
the following conditions:
\begin{itemize}
\item[i)] Each $A_j$ is a fuzzy set in $X$ and
$\sup_j\{\mu_{A_j}(x)\}=1$, for all $x\in X$.
\item[ii)] Each $\phi_J$ is a
bijection, defined on the support of $A_j$,
$$\{x\in X: \mu_{A_j}(x)>0\},$$
which  maps $A_j$ onto an open fuzzy set $\phi_j[A_j]$ in some
fuzzy topological vector space $E_j$, and, for each $l\in J$,
$\phi_j[A_j\cap A_l]$ is an open fuzzy set in $E_j$.
\item[iii)]  The mapping $\phi_l\circ\phi_j^{-1}$, which maps $\phi_j[A_j\cap
A_l]$ is a ${\cal C}^1-$fuzzy diffeomorphism for each pair of
indices $j,l$.
\end{itemize}
Each pair $(A_j,\phi_j)_{j\in J}$ is called {\it a fuzzy chart} of
the fuzzy atlas. If a point $x\in X$ lies in the support of $A_j$
then $(A_j,\phi_j)_{j\in J}$ is said to be a fuzzy chart at $x$.

Let $(X,\tau)$ be a fuzzy topological space. Suppose there exists
an open fuzzy set $A$ in $X$ and a fuzzy continuous bijective
mapping $\phi$ defined on the support of $A$ and mapping onto an
open fuzzy set $V$ in some fuzzy topological vector space $E$.
Then $(A,\phi)$ is said to be compatible with the ${\cal
C}^1-$atlas $\{(A_j,\phi_j)\}_{j\in J}$ if each mapping
$\phi_j\circ\phi^{-1}$ of $\phi[A\cap A_j]$ onto $\phi_j[A\cap
A_j]$ is a ${\cal C}^1-$fuzzy diffeomorphism. Two ${\cal
C}^1-$fuzzy atlases are compatible if each fuzzy chart of one
atlas is compatible with each fuzzy chart of the other atlas. It
may be verified immediately that the relation of compatibility
between ${\cal C}^1-$fuzzy atlases is an equivalence relation. An
equivalence class of $C^{1}$ fuzzy atlases on $X$ is said to
define a ${\cal C}^1-$fuzzy manifold on $X$.
\paragraph{Proposition 3.3.}
Let $X,Y$ be the fuzzy manifolds; then the product $X\times Y$ is
a fuzzy manifold.
\paragraph{Definition 3.4.}
Let $X,Y$ be the fuzzy manifolds and let $f$ be a mapping  of $X$
into $Y$. Then $f$ is said to be fuzzy differentiable at a point
$x\in X$ if there is a fuzzy chart $(U,\phi)$ at $x\in X$ and a
fuzzy chart $(V,\varphi)$ at $f(x)\in Y$ such that the mapping
$\varphi\circ f\circ\phi^{-1}$, which maps $\phi[U\cap f^{-1}[V]]$
into $\varphi[V]$ is fuzzy differentiable at $\phi(x)$. The
mapping $f$ is fuzzy differentiable if it is fuzzy differentiable
at every point of $X$; it is a ${\cal C}^1-$fuzzy diffeomorphism
if $\varphi\circ f\circ\phi^{-1}$ is a ${\cal C}^1-$fuzzy
diffeomorphism.
\paragraph{Example 3.5.}
Suppose that $X={\SS}^1$ is the set of points of the unit circle
in ${\RR}^2$. If $U$ is the fuzzy set of ${\SS}^1$ consisting of
the points
$$(\sin2\pi t,\cos2\pi t),~~~~~0<s<1,$$
with the characteristic function $\mu_U:{\SS}^1\longrightarrow I$
defined by
$$\mu_U(\sin2\pi t,\cos2\pi t)=\sin^22\pi t+\cos^22\pi t=1,$$
then the function
\begin{equation}
{\phi_1:{\SS}^1\longrightarrow {\RR}\atop (\sin2\pi t,\cos2\pi
t)\longmapsto t} \hskip 1cm
\end{equation}
is a bijection onto an open fuzzy set of ${\RR}$ and so
$(U,\phi_1)$ is a fuzzy chart for ${\SS}^1$. If $V$ be another
fuzzy set
$$(\sin2\pi t,\cos2\pi t),~~~~~-{1\over2}<t<{1\over2},$$
of ${\SS}^1$ with the characteristic function
$\mu_V:{\SS}^1\longrightarrow I$ defined by
$$\mu_U(\sin2\pi t,\cos2\pi t)={1\over2},$$
the function
\begin{equation}
{\phi_2:{\SS}^1\longrightarrow {\RR}\atop (\sin2\pi t,\cos2\pi
t)\longmapsto t} \hskip 1cm
\end{equation}
is another such chart. Each of $U$ and $V$ are fuzzy sets of
${\SS}^1$ such that
$$\sup\{\mu_U(x),\mu_V(x)\}=1.$$
Since
$$\phi_1=\phi_2~~if~~0<\phi_1<{1\over2}~~~~and~~~~\phi_2=\phi_1-1~~if~~{1\over2}<\phi_1<1,$$
The mapping $\phi_2\circ\phi_1^{-1}:{\RR}\longrightarrow{\RR}$ is
defined by
$$\phi_2\circ\phi_1^{-1}(t)=\left\{ {t \atop t-1}     \hskip 1cm
{0<t<{1\over2}, \atop {1\over2}<t<1.} \right.$$

Clearly this mapping is a ${\cal C}^1-$fuzzy diffeomorphism,
therefor the fuzzy charts $U$ and $V$ form a ${\cal C}^1-$fuzzy
atlas on ${\SS}^1$. Now we define another ${\cal C}^1-$fuzzy atlas
on ${\SS}^1$. Let $U_1$ be the set of points
$$(z_1,z_2)\in{\SS}^1~~~,~~~z_1>0.$$
We have $U_1$ a fuzzy set of ${\SS}^1$ by the characteristic
function $\mu_{U_1}:{\SS}^1\longrightarrow I$ defined by
$$\mu_{U_1}(z_1,z_2)={1\over4}.$$
Also the function
\begin{equation}
{\psi_{U_1}:{\SS}^1\longrightarrow {\RR}\atop (z_1
,z_2)\longmapsto z_2,} \hskip 1cm
\end{equation}
is a bijection onto an open fuzzy set of ${\RR}$. It is therefore
a fuzzy chart for ${\SS}^1$. Let $U_2,U_3,U_4$ be the sets of
points of ${\SS}^1$ such that $z_2>0, z_1<0, z_2>0$ respectively.
They are fuzzy set by the characteristic functions
$\mu_{U_i}:{\SS}^1\longrightarrow I$ that $i=2,3,4$ by
$$\mu_{U_2}(z_1,z_2)={1\over4}~,~\mu_{U_3}(z_1,z_2)={1\over4}~,~\mu_{U_4}(z_1,z_2)={1\over4}.$$
By the following bijection onto an open fuzzy set of ${\RR}$,
$$\psi_{U_2}(z_1,z_2)=z_1~,~\psi_{U_3}(z_1,z_2)=z_2~,~\psi_{U_4}(z_1,z_2)=z_1$$
$U_1,U_2,U_3,U_4$ are the fuzzy charts which cover ${\SS}^1$ and
for any intersections the change of coordinates is a ${\cal C}^1$
fuzzy diffeomorphism. For instance,
$$\psi_2=\sqrt{1-\psi_1^2}~~~,~~~\psi_1>0, $$
and so the mapping
\begin{equation}
\psi_2\circ\psi_1^{-1}:I\longrightarrow I\atop t\longmapsto
\sqrt{1-t^2}, \hskip 1cm
\end{equation}
is a ${\cal C}^1-$fuzzy diffeomorphism on the open interval
$(0,1)$. Since, we have
$$\psi_2=-\sqrt{1-\psi_3^2}~~~,~~~\psi_3>0, $$
$$\psi_3=\sqrt{1-\psi_1^2}~~~,~~~\psi_1<0, $$
$$\psi_4=\sqrt{1-\psi_3^2}~~~,~~~\psi_3<0, $$
then the other changes coordinates are also ${\cal C}^1-$fuzzy
diffeomorphisms. Therefore the charts
$\psi_1,\psi_2,\psi_3,\psi_4$ form a ${\cal C}^1-$fuzzy atlas. In
fact, this ${\cal C}^1-$fuzzy atlas is compatible with the
pervious ${\cal C}^1-$fuzzy atlas. To prove this, we need to show
that, the additional change of coordinates are also ${\cal C}^1-$
fuzzy diffeomorphism. This following since

$$\psi_1=\left\{ {\cos(2\pi\phi_1) \atop \cos(2\pi\phi_2)}     \hskip 1cm
{0<\phi_1<{1\over2}, \atop 0<\phi_2<{1\over2},} \right.$$

$$\psi_2=\left\{ {\cos(2\pi\phi_1) \atop \cos(2\pi\phi_2)}     \hskip 1cm
{0<\phi_1<{1\over4}~{\rm or}~{3\over4}<\phi_1<1, \atop
-{1\over4}<\phi_2<{1\over4},} \right.$$
 and so on.

\paragraph{Definition 3.6.}
A fuzzy differentiable function $\psi:M^{\prime}\longrightarrow M$
is called a {\it fuzzy immersion} if its rank is equal to the
dimension of $M^{\prime}$ at each point of its domain. If its
domain is the whole of $M^{\prime}$, $\psi$ is said to be a fuzzy
immersion of $M^{\prime}$ into $M$.

\paragraph{Definition 3.7.} A ${\cal C}^1-$fuzzy manifold $M^{\prime}$ is
said to be a ${\cal C}^1-$fuzzy submanifold of a ${\cal C}^1-$
fuzzy manifold $M$ if

(i) $M^{\prime}$ is a fuzzy subset of $M$,

(ii) Natural fuzzy injection $j:M^{\prime}\longrightarrow M$ is a
fuzzy immersion.

\paragraph{Example 3.8.} Clearly $M(n\times n,{\RR})$, the set
of real $n\times n$ matrices, is a ${\cal C}^1-$fuzzy manifold and
${\rm GL}(n,{\RR})$ ia a fuzzy subset of it. If $j: {\rm
GL}(n,{\RR})\longrightarrow M(n\times n,{\RR})$ is the natural
fuzzy injection and
 $$\det\circ j: {\rm GL}(n,{\RR})\longrightarrow{\RR}$$ is fuzzy differentiable, then $j$ is a
fuzzy immersion so $ {\rm GL}(n,{\RR})$ is a ${\cal C}^1-$fuzzy
submanifold of $M(n\times n,{\RR})$.
\section{ Fuzzy  Lie group }
\paragraph{Definition 4.1.}
{\it A fuzzy Lie group $G$} is a ${\cal C}^1-$fuzzy manifold $G$
which is also $G$ is a group, such that the mappings
$$m:(G\times G,\tau\times\tau)\longrightarrow (G,\tau),\;\;\;\;i: (G,\tau)\longrightarrow (G,\tau)$$
defined in (\ref{eq:1}) and (\ref{eq:2}), be fuzzy differentiable.
\paragraph{Example 4.2.}
(i) One of the simplest example of a Lie group fuzzy is ${\RR}^n$
that is commutative fuzzy Lie group. The group operation is given
by vector addition . The identity element is the zero vector, and
the inverse of a vector $x$ is the vector $-x$. If ${\RR}^n$
equipped with the ordinary fuzzy topology, it is trivial which the
mappings
\begin{eqnarray*}
&m:{\RR}^n\times{\RR}^n\longrightarrow{\RR}^n \hspace{1cm}
(x,y)\mapsto x+y,& \\ &
i:{\RR}^n\longrightarrow{\RR}^n\hspace{1cm} x\mapsto x^{-1},&
\end{eqnarray*}
are fuzzy differentiable.

(ii) The other example of fuzzy lie group is the general linear
group ${\rm GL}(n,{\RR})$ consisting of all invertible $n\times n$
real matrices, with matrix multiplication defining the group
multiplication, and matrix inversion defining the inverse. In
fact, ${\rm GL}(n,{\RR})$ is an $n^2$-dimensional ${\cal
C}^1-$fuzzy manifold such that
\begin{eqnarray*}
&{m:{{\rm GL}(n,{\RR}^n)\times {\rm GL}(n,{\RR}^n)\longrightarrow
{\rm GL} (n,{\RR}^n)}\hspace{1cm} {(A,B)\mapsto AB},}&
\\ &{i:{{\rm GL}(n,{\RR}^n)\longrightarrow
 {\rm GL}(n,{\RR}^n)}\hspace{1cm}  {A\mapsto A^{-1}},}&
\end{eqnarray*}
\paragraph{Definition 4.3.}
A fuzzy transformation group acting on a ${\cal C}^1-$fuzzy
manifold is determined by a fuzzy Lie group $G$ and fuzzy
differentiable map $\Phi:G\times M\longrightarrow M$, which
satisfies
\begin{itemize}
\item[i)] $\Phi$ is a fuzzy global surjective,
\item[ii)] $\Phi(g,\Phi(h,x))=\Phi(gh,x),$ for any $x\in M$ and
$g,h\in G$.
\end{itemize}
\paragraph{Example 4.4.}
${\rm GL}(n,{\RR}^n)$ acts on ${\RR}^n$ as a fuzzy transformation
group with the map
$$\Phi:{{\rm GL}(n,{\RR}^n)\times{\RR}^n\longrightarrow
{\RR}^n}\hspace{1cm} (A,x)\mapsto Ax.$$
\paragraph{Definition 4.5.}
Let $G$ be a fuzzy Lie group, then $H\subset G$ is called {\it a
fuzzy Lie subgroup} if $H$ is both a subgroup and a ${\cal
C}^1-$fuzzy submanifold. For instance, $O(n,{\RR})$, that is real
orthogonal $n\times n$ matrices, is a fuzzy Lie subgroup of ${\rm
GL} (n,{\RR})$.
\paragraph{Proposition 4.6.}
If a fuzzy Lie group $G$ acts on ${\cal C}^1-$fuzzy manifold $M$
as a fuzzy transformation group then so does any fuzzy Lie
subgroup $H$ of $G$.

\medskip \noindent {\it Proof:} If $j:H\longrightarrow G$ is the natural
fuzzy injection and $id: M\longrightarrow M$ be a fuzzy identity
map. There is a suitable global function $\Phi_H:H\times
M\longrightarrow M$ such that $$\Phi_{H}=\Phi\circ(j\times i),$$
where $j\times i:H\times M\longrightarrow G\times M.$ Therefore
$\Phi_H$ is a fuzzy surjective and fuzzy differentiable function
which we have
\begin{eqnarray*}
{\Phi_H}(h,{\Phi_H}(h^{\prime},x)) &=& ({\Phi}\circ{(j\times
i)})\bigg(h,({\Phi}\circ{(j\times i))(h^{\prime},x)}\bigg)\\
&=&\Phi(h,\Phi(h^{\prime},x)) =\Phi(hh^{\prime},x)\\
&=&(\Phi\circ(j\times i))(hh^{\prime},x) = {\Phi_H}(hh^{\prime},x)
\end{eqnarray*}
for any $x\in M$ and $h,h^{\prime}\in H$.\hfill\
\paragraph{Example 4.7.}
In the example 3.4., we se that the ${\rm GL}(n,{\RR})$ be a fuzzy
Lie group, and $O(n,{\RR})$ is one of its fuzzy Lie subgroups,
then  $O(n,{\RR})$ also acts on ${\RR}^3$ as a fuzzy
transformation group.
\paragraph{Definition 4.8.}
Let $\Phi:G\times M\longrightarrow M$ be a fuzzy transformation
group. A subset $S\subset M$ is called {\it fuzzy $G$-invariant
subset of $M$}, if $\Phi(G\times S)\subseteq S$.
\paragraph{Proposition 4.9.}
Consider $\Phi:G\times M\longrightarrow M$ as a fuzzy
transformation group. If a regular ${\cal C}^1-$fuzzy submanifold
$M^{\prime}$ of ${\cal C}^1-$fuzzy manifold $M$ is $G$-invariant,
then $G$ acts naturally on $M^{\prime}$ as a fuzzy transformation
group.

\medskip \noindent {\it Proof:} Suppose that $k:G\times M^{\prime}\longrightarrow
G\times M$ is the natural fuzzy injection, then the function
$\Phi^{\prime}:G\times M^{\prime}\longrightarrow M^{\prime}$
induced by $\Phi\circ k$ is fuzzy differentiable that is defines
the required action of $G$ on $M^{\prime}$.\hfill\
\paragraph{Proposition 4.10.}
If ${M/\rho}$ is a quotient ${\cal C}^1-$fuzzy manifold of $M$ and
equivalence relation $\rho$ is preserved by a fuzzy Lie
transformation group $G$ on $M$ then $G$ acts naturally on
${M/\rho}$.

\medskip \noindent {\it Proof:} Let $\alpha:M\longrightarrow{M/\rho} $ be
natural fuzzy surjection and $\Phi:G\times M\longrightarrow M$ be
fuzzy differentiable map so $\alpha\circ\Phi:G\times
M\longrightarrow {M/\rho}$ is a fuzzy differentiable function.
$G\times({M/\rho})$ is a quotient ${\cal C}^1-$fuzzy manifold of
$G\times M$ and $\alpha\circ\Phi$ is an invariant of corresponding
equivalence relation on $G\times M$. It therefore projects to the
fuzzy differentiable function
$${\Phi:G\times{M/\rho}\longrightarrow
{M/\rho}\hspace{1cm} {(g,\alpha m)\mapsto \alpha(gm)}.}$$ This
defines the required action $G$ on ${M/\rho}$.\hfill\

\end{document}